# Operators in the mind: Jan Lukasiewicz and Polish notation

Eduardo Mizraji
Group of Cognitive Systems Modeling, Biophysics and Systems Biology Section
Facultad de Ciencias, Universidad de la República,
Montevideo, Uruguay
Correspondence: emizraji@gmail.com; mizraj@fcien.edu.uy
ORCID id - 0000-0001-6938-8427

ABSTRACT
In 1929 Jan Lukasiewicz used, apparently for the first time, his Polish notation to represent the operations of formal logic. This is a parenthesis-free notation, which also implies that logical functions are operators preceding the variables on which they act. In the 1980s, within the framework of research into mathematical models on the parallel processing of neural systems, a group of operators emerged -neurally inspired and based on matrix algebra- which computed logical operations automatically. These matrix operators reproduce the order of operators and variables of Polish notation. These logical matrices can also generate a three-valued logic with broad similarities to Lukasiewicz's three-valued logic. In this paper, a parallel is drawn between relevant formulas represented in Polish notation, and their counterparts in terms of neurally based matrix operators. Lukasiewicz's three-valued logic, shown in Polish notation has several points of contact with what matrices produce when they process uncertain truth vectors. This formal parallelism opens up scientific and philosophical perspectives that deserve to be further explored.

## 1. Introduction

Jan Lukasiewicz was a man of great creativity and a rare insight. In 1963, the English translation of his book "Elements of Mathematical Logic" (Lukasiewicz, 1929) was published, having first appeared in Polish in 1929. This book arose from the notes of his lectures delivered at Warsaw University in the academic year 1928-1929. In the preface to the original edition, Lukasiewicz makes a list of what he believed (with modesty) to be his contributions to the theory of logic. As his first choice, he states: "*1. The parenthesis-free notation of expressions in the sentential calculus and in Aristotle´s syllogistic*". Then he mentions the axioms of sentential calculus and in fifth place the systems of many-valued logic (Lukasiewicz 1929 [1963, p. IX]).

Both in his article "Creative Elements in Science" of 1912 and in his "Farewell Lecture" of 1918, Lukasiewicz vindicated the creative freedom of the researcher in logic. In "Creative Elements" he writes: "*The a priori mental constructions, which are contained in every synthesis, imbue the whole science ideal and creative elements*" (Lukasiewicz 1912). In his "Farewell Lecture" he said: "*The possibility of constructing different logical systems shows that logic is not restricted to reproduction of facts but is a free product of man, like a work of art*" (Lukasiewicz 1918).

Despite these opinions, Lukasiewicz was a convinced opponent of the idea that formal logic was a branch of psychology, and he expressed this anti-psychologism position in logic with emphasis in his 1929 text.

The aim of this paper is first to describe the parenthesis-free notation that Lukasiewicz prioritized in enumerating his achievements, then to point out the intense originality of that invention, and after that to show a remarkable situation: The strong (and perhaps unexpected) connections of Polish notation with a mathematical format that emerged from the territory of biological memories.

We think Lukasiewicz was probably right that formal logic is not a branch of psychology. But our time perhaps shows that formal logic can be a corollary of neural computation.

## 2. Polish notation

His 1929 text (as far as we know) seems to be the first place where Lukasiewicz sets out his Polish notation. He does so without giving any indication of how he arrived at it. In Chapter 2 of his book ("The sentential calculus") he describes it with an elementary arithmetical example, and in the process shows its generality. He writes
a + (b + c) = (a + b) + c
and shows that this expression is equivalent to
+ a + bc = + + abc.
Thus, the symbol "+" is a (dyadic) operator that acts on the operators and expressions located to its right (Lukasiewicz 1928 [1963, p. 24]). If there is a second operator to the right of the first operator, the first variable on which the first operator acts is the one resulting from the computation performed by the second.

With this formalism, he defines the basic logical operations, negation N, implication C, conjunction K and disjunction A. He defines the truth values as 0 for "false" and 1 for "true". Their definitions are:
N0 = 1 ; N1 = 0,
C11 = C01 = C00 =1; C10 =0,
K11 = 1; K10 = K01 = K00 = 0.
A11 = A10 = A01 =1; A00 =0.
The classical definition of implication in terms of disjunction and negation is written as
Cpq = ANpq.
Variables represented by lowercase letters belong to the set $\{0, 1\}$.
De Morgan's laws in Polish notation have the following format:
Kpq = NANpNq,
Apq = NKNpNq.
We will now show some basic laws of formal logic using Polish notation:
*Principle of Non-Contradiction*: NKpNp
*Principle of the excluded middle:* ApNp
*Reduction to absurdity:* CCpKqNqNp
*Hypothetical syllogism:* CKCpqpq

In Bochenski (1959) numerous fundamental formulas of logical calculus are shown both in Polish notation and in Whitehead-Russell notation.
We conclude this section by pointing out that in his 1929 text Lukasiewicz develops a singular argument to justify the structure of implication. The argument goes in this direction: If x is divisible by 9, then x is divisible by 3. If this implication is assumed to be true for all x, then implications are true where x = 18 (C11 = 1), also if x = 15 (C01 = 1) and also if x = 19 (C00 = 1). He then adds that if it is not possible to find a substitution that makes the antecedent true and the consequent false then an implication of that kind is not true (C10 = 0). This last argument and its apparent circularity (of using an implicit implication to define a property of implication) should possibly be adopted as an axiom.

Presumably a similar approach influenced Lukasiewicz's definition of trivalent implication and led to non-coincidence with other definitions of implication when both truth values were uncertain (e.g. in the trivalent logic of Kleene, 1938).

## 3. The foundations of neural computing

The origin of the current theory of neural computing, after many important antecedents, was consolidated in the early 1970s and is based on the concept of "neural vector". This concept arises from the fact that in complex nervous systems, such as the human one, neural information is encoded by means of an extensive set of signals that are processed in parallel by various regions of the brain. This extensive set of parallel signals is adequately represented by mathematical vectors.

In 1972, two independent papers were published that converged on the same formalism. The authors' goal was to show how information encoded using vectors was stored in memories. These memories were represented by matrices and showed remarkable properties, such as distributed storage of data and partial resistance to deterioration of the matrix, among others (Anderson 1972; Kohonen 1972). In these models, the matrix coefficients represented the synaptic connections between neurons, and these models suggested that these synapses could be the physical residence of the information stored in memory. Another fundamental aspect of these models was that they were capable of being trained in such a way as to incorporate new associations into memory. An exploration of this matrix formalism showed how concepts could be interpreted as average vectors, representing prototypes generated by a set of related items (Cooper 1973).

In 1982, Hopfield published an article where memories are such that vectors "search" for their associations through a dynamic process (Hopfield 1982).
It soon became apparent that these models had severe limitations in modulating their associations by context, a natural and necessary modulation in biological memories. An elegant mathematical proof of these limitations is given in Hinton (1989).
Given these limitations of the aforementioned models, in the middle of that decade articles were published that introduced complexities into the models of previous years, such as hidden layers of neurons, which extended the learning and recognition capabilities of previous memory models (see, for example, Rumelhart, Hinton and McClelland 1986).

On the other hand, seeking to overcome the problem of contexts but retaining the mathematics of matrix models (and their potential to develop further theories), tensor operations were introduced between the neural vectors that acted as basic patterns and other neural vectors that acted as contexts. This tensor representation was described, almost coincidentally in time, by several independent researchers (Humphreys, Bain and Pike 1989; Mizraji 1989; Smolensky 1990).

It should be noted that, in recent years, these investigations on neural models with learning capacity from the 1980s "exploded" into the powerful and varied models of Artificial Intelligence. These AI models are fully based on those findings but take them to another level of complexity allowed by the enormous computational capacities that exist today (Strink and Kondratenko 2021; Valle-Lisboa, Pomi and Mizraji 2023).

In Mizraji (1989) a case was described (the exclusive-or, an important function for exploring the contextualization capabilities in neuronal memories) that shows that this context-dependent matrix memory formalism allowed the representation of logical operations. On this basis, several papers were published, showing that every well-formed formula of basic logic could be represented by a formalism based on matrix operators and vectors (see, for example, Mizraji 1992, 1996, 2008).

In the next section, we will show how Lukasiewicz's Polish notation has a striking parallel to the order of matrix and vector operators that arises from the above-mentioned formalism based on context-dependent matrix memories. This leads us to emphasize that Lukasiewicz's Polish notation is primarily an operator theory and that the fact that it is parenthesis-free is a consequence of this primordial fact.

## 4. Polish notation and "neural logic"

Our notation differs from that used by Mizraji (1992) because here we will adopt the Polish symbols for the operators. We will assume that the value 0 used by Lukasiewicz for "false" corresponds to a vector $\mathbf{f}$, and that the value 1, "true", corresponds to the vector $\mathbf{t}$. These new truth values are column vectors of dimension 2 or higher, and we will assume that the vectors $\mathbf{f}$ and $\mathbf{t}$ are orthogonal to each other, and of unit module (orthonormal). These assumptions facilitate calculations because scalar products, extremely important in this formalism, are simplified for orthonormal vectors. Thus, the scalar products between truth vectors are $(\mathbf{f}.\mathbf{f}) = (\mathbf{t}.\mathbf{t}) = 1$ and $(\mathbf{f}.\mathbf{t}) = (\mathbf{t}.\mathbf{f}) = 0$. In addition to scalar products, the operation of the matrix version of logic uses matrix products (which include matrix-vector products) and tensor products. Details about logical matrices and the products mentioned are given in the Appendix.

Logical matrices are actually matrix memories that store the properties of basic logical operations. Here we will use for the matrix operators the same notation (but in bold) that Lukasiewicz used for his logical functions Negation, Implication, Conjunction and Disjunction: $\mathbf{N}$, $\mathbf{C}$, $\mathbf{K}$ and $\mathbf{A}$, respectively. The formulas that we will show below must be written using the parentheses hierarchy required by matrix calculus, but the similarity with Polish notation is deeper, as we will see, than the parenthesis-free nature, since the matrix formalism gives the interpretation of Polish notation another dimension. The symbol $\otimes$ denotes the tensor products connecting the vector truth values.

We now show how matrices operate on truth vectors $\mathbf{f}$ and $\mathbf{t}$ and display the various matrix logic formulas following the same order we used in Section 2 to illustrate Polish notation:

$\mathbf{Nf} = \mathbf{t}$ ; $\mathbf{Nt} = \mathbf{f}$,

$\mathbf{C}(\mathbf{t} \otimes \mathbf{t}) = \mathbf{C}(\mathbf{f} \otimes \mathbf{t}) = \mathbf{C}(\mathbf{f} \otimes \mathbf{f}) = \mathbf{t}$; $\mathbf{C}(\mathbf{t} \otimes \mathbf{f}) = \mathbf{f}$,

$\mathbf{K}(\mathbf{t} \otimes \mathbf{t}) = \mathbf{t}$; $\mathbf{K}(\mathbf{t} \otimes \mathbf{f}) = \mathbf{K}(\mathbf{f} \otimes \mathbf{t}) = \mathbf{K}(\mathbf{f} \otimes \mathbf{f}) = \mathbf{f}$.

$\mathbf{A}(\mathbf{t} \otimes \mathbf{t}) = \mathbf{A}(\mathbf{f} \otimes \mathbf{t}) = \mathbf{A}(\mathbf{t} \otimes \mathbf{f}) = \mathbf{t}$; $\mathbf{A}(\mathbf{f} \otimes \mathbf{f}) = \mathbf{f}$.

In what follows we use lowercase bold letters to designate vectors in the set $\{\mathbf{f}, \mathbf{t}\}$. Implication $\mathbf{C}$ complies with its classical formula

$\mathbf{C} = \mathbf{A}(\mathbf{Nu} \otimes \mathbf{v})$.

The matrix De Morgan laws are:

$\mathbf{K}(\mathbf{u} \otimes \mathbf{v}) = \mathbf{NA}(\mathbf{Nu} \otimes \mathbf{Nv})$,

$\mathbf{A}(\mathbf{u} \otimes \mathbf{v}) = \mathbf{NK}(\mathbf{Nu} \otimes \mathbf{Nv})$.

The basic laws of formal logic shown above take the following formats in the matrix framework:

*Principle of Non-Contradiction*: $\mathbf{NK}(\mathbf{p} \otimes \mathbf{Np})$

*Principle of the excluded middle*: $\mathbf{A}(\mathbf{p} \otimes \mathbf{Np})$

*Reduction to absurdity*: $\mathbf{C}\{\mathbf{C}[\mathbf{p} \otimes \mathbf{K}(\mathbf{q} \otimes \mathbf{Nq}] \otimes \mathbf{Np}\}$

*Hypothetical syllogism*: $\mathbf{C}\{\mathbf{K}[\mathbf{C}(\mathbf{p} \otimes \mathbf{q}) \otimes \mathbf{p}] \otimes \mathbf{q}\}$

We see that if we ignore the parentheses, the (implicit) matrix multiplications and the tensor product signs, the order of operators and variables is the same in these last formulas and in those written in Section 2 with Polish notation. Let us also note that this remarkable fact is an inevitable consequence of the use of matrix formalism since it is imposed by the operating rules of linear algebra.

This coincidence in the order of operators and variables is very surprising when we consider that Lukasiewicz created his Polish notation when the use of matrices was limited to a few domains of science (even if Peirce had made a precursory use of matrices for his theory of relations (cited in Copilowich 1948)). And, obviously, the Polish notation preceded by many decades the development of neural computation models. Already in 1929 Lukasiewicz had a clear perception of the power, and perhaps the naturalness, of the use of operators to formalize logical calculus, an operator formalism from which he did not detach himself in the rest of his investigations.
But the analysis of the three-valued logic created by Lukasiewicz and the way logical matrices process uncertain truth values, holds more surprises.

## 5. Three-valued logic and matrix operators
Lukasiewicz describes his three-valued logic in a short paper published in 1920, where he does not use Polish notation (Lukasiewicz 1920). The intellectual mechanisms that led Lukasiewicz to assign truth values to logical functions when faced with uncertain truth values do not seem to have been made explicit (Urquhart 2021). Nor is how Kleene (1938) created his own version of his trivalent logic.

Given the variety of mathematical talents, one can conjecture the use of two different methods for assigning truth values to uncertain variables. On the one hand, it is possible to reason in stages. Let us look at a simple case. Given the disjunction **A**, if a truth value is to be assigned to **A**(1,?) where it is unknown whether ? is 0 or 1, then one can be certain that in any case **A**(1,?) = 1. But if **A**(0,?), depending on whether the uncertainty is 1 or 0, the result diverges and becomes uncertain, so that **A**(0,?) = ?. But on the other hand, people with extraordinary mathematical abilities "see" the results without being (at least consciously) analytical. In this regard, it should be noted that Hadamard deeply investigated the intellectual styles of mathematical creation (Hadamard 1945).

In his three-valued logic, Lukasiewicz symbolized uncertain truth values using the symbol ½. Let us now show the evaluations that must be added to the classical operations to define three-valued operators.

*1) Negation*
N ½ = ½ .
*2) Implication*
C 1 ½ = ½ , C ½ 1 = 1,
C 0 ½ = 1 , C ½ 0 = ½ ,
C ½ ½ = 1.
*3) Conjunction*
K 1 ½ = ½ , K ½ 1 = ½ ,
K 0 ½ = 0 , K ½ 0 = 0 ,
K ½ ½ = ½ .
*4) Disjunction*
A 1 ½ = 1 , A ½ 1 = 1,
A 0 ½ = ½ , A ½ 0 = ½ ,
A ½ ½ = ½ .

To evaluate how matrices handle logical uncertainties, let us define a vector **i** where the truth values **f** and **t** are both weighted by the coefficient 1/2:

$$\mathbf{i} = (1/2)\mathbf{f} + (1/2)\mathbf{t}$$

Since matrices are linear operators, the weighted vectors are processed naturally. Here we will directly show the results associated with the evaluations shown above for the Lukasiewicz operators. The Appendix shows a general method described in Mizraji (1992) that allows immediate calculations for any weighting of the vectors $\mathbf{f}$ and $\mathbf{t}$.

*1) Negation*
$$\mathbf{Ni} = \mathbf{i} \; .$$

*2) Implication*
$$\mathbf{C}(\mathbf{t} \otimes \mathbf{i}) = \mathbf{i} \;, \quad \mathbf{C}(\mathbf{i} \otimes \mathbf{t}) = \mathbf{t} \;,$$
$$\mathbf{C}(\mathbf{f} \otimes \mathbf{i}) = \mathbf{t} \;, \quad \mathbf{C}(\mathbf{i} \otimes \mathbf{f}) = \mathbf{i}$$
$$\mathbf{C}(\mathbf{i} \otimes \mathbf{i}) = (1/4)\mathbf{f} + (3/4)\mathbf{t} \; .$$

*3) Conjunction*
$$\mathbf{K}(\mathbf{t} \otimes \mathbf{i}) = \mathbf{i} \;, \quad \mathbf{K}(\mathbf{i} \otimes \mathbf{t}) = \mathbf{i} \;,$$
$$\mathbf{K}(\mathbf{f} \otimes \mathbf{i}) = \mathbf{f} \;, \quad \mathbf{K}(\mathbf{i} \otimes \mathbf{f}) = \mathbf{f}$$
$$\mathbf{K}(\mathbf{i} \otimes \mathbf{i}) = (3/4)\mathbf{f} + (1/4)\mathbf{t}$$

4) *Disjunction*
$$\mathbf{A}(\mathbf{t} \otimes \mathbf{i}) = \mathbf{t} \;, \quad \mathbf{A}(\mathbf{i} \otimes \mathbf{t}) = \mathbf{t} \;,$$
$$\mathbf{A}(\mathbf{f} \otimes \mathbf{i}) = \mathbf{i} \;, \quad \mathbf{A}(\mathbf{i} \otimes \mathbf{f}) = \mathbf{i}$$
$$\mathbf{A}(\mathbf{i} \otimes \mathbf{i}) = (3/4)\mathbf{f} + (1/4)\mathbf{t}$$

We see that the negation matrix $\mathbf{N}$ operates in a manner analogous to the three-valued negation. The matrix outcomes of $\mathbf{C, K}$ and $\mathbf{A}$ coincide in their vector version for all pairs of variables except for the pair ($\mathbf{i,i}$), where linear combinations with different biases arise. It is interesting that Kleene's logic for two uncertain values yields one uncertain value (Kleene 1938). Thus, the curious fact arises that the weighting of the vector $\mathbf{t}$ in the matrix evaluation of the implication, is an average of the evaluation 1 of Lukasiewicz and the ½ of Kleene.

But if we focus on the convergent aspects between Lukasiewicz's three-valued logic and matrix computations, except for the value of the double uncertainty, the overlap is complete, which can be surprising.

## 6. Concluding remarks

By inventing Polish notation, Lukasiewicz created a logic of operators that many decades later was structurally reproduced in the field of neural computing. Polish notation does not use parentheses, and the link between operators and variables is sufficient to define operations unambiguously. However, the notation of logical operations arising from neural models, which exactly reproduces the order of operators and variables of Polish notation, does require parentheses and signs (such as the tensor product) in accordance with the uses of matrix algebra.

If it were agreed that operators are matrices and that adjacent vectors are connected by tensor products, a parenthesis-free matrix-vector representation would be recovered. However, the powerful factorization capabilities of the Kronecker product (not shown here, but discussed in Mizraji 1996) would be difficult without tensor notation.

But there is another deeper argument to establish the link between Polish notation and the neural representation of logical functions. The algebraic representation that aims to imitate neurocomputational processes uses parentheses and signs on the surface of paper or a screen, that is, *in its symbolic writing*. If processes of this type, perhaps more complex but summarizable through the neural models mentioned, really occur in the mind, there the computations do not operate through symbolism that requires the spelling of mathematics. Matrix or tensor products arise from an analog computation process physically guided by the interactions of neural signals with the synaptic junctions of memory modules.

It is very likely that human beings have these operators in their minds. The neural structures that support the memories on which the biographical identity of each person is based, the syntactic structure of language and also the semantic networks, are a strong argument in favor of our mind being structured on memory modules. And it is plausible that (congenital or acquired, and without entering into the "Platonic-Aristotelian" controversy) an adult human has incorporated logical modules that are activated under conditions of formal reasoning (e.g.: during the demonstration of elementary theorems in the early stages of education, which are usually learned before formal logic).

The fact that, at the same time, the matrix format mimics the order of operators and variables of Polish notation and generates a trivalent logic, so close to that created by Lukasiewicz, is interesting and curious. If these matrix memory models touch on aspects of neural reality, then the coincidence described between Polish notation and matrix formalism would indicate an extraordinary intuition of Jan Lukasiewicz about the modes of operation of the human mind.

## Appendix

The details of the operations that we will show in this Appendix are in matrix algebra texts such as Barnett (1990).

The matrices associated with the logical operations that we mention in the text are the following:

$\mathbf{N} = \mathbf{f}\mathbf{t}^T + \mathbf{t}\mathbf{f}^T$,

$\mathbf{C} = \mathbf{t}(\mathbf{t} \otimes \mathbf{t})^T + \mathbf{f}(\mathbf{t} \otimes \mathbf{f})^T + \mathbf{t}(\mathbf{f} \otimes \mathbf{t})^T + \mathbf{t}(\mathbf{f} \otimes \mathbf{f})^T$,

$\mathbf{K} = \mathbf{t}(\mathbf{t} \otimes \mathbf{t})^T + \mathbf{f}(\mathbf{t} \otimes \mathbf{f})^T + \mathbf{f}(\mathbf{f} \otimes \mathbf{t})^T + \mathbf{f}(\mathbf{f} \otimes \mathbf{f})^T$,

$\mathbf{A} = \mathbf{t}(\mathbf{t} \otimes \mathbf{t})^T + \mathbf{t}(\mathbf{t} \otimes \mathbf{f})^T + \mathbf{t}(\mathbf{f} \otimes \mathbf{t})^T + \mathbf{f}(\mathbf{f} \otimes \mathbf{f})^T$.

Vectors $\mathbf{f}$ and $\mathbf{t}$ are column vectors with dimension $d \geq 2$. The superscript $\mathbf{T}$ indicates transpose, which converts column vectors into rows. The matrix $\mathbf{N}$ is a square matrix of dimension $d \times d$; dyadic matrices are rectangular of dimension $d \times d^2$.

The tensor product $\otimes$ (or Kronecker product) has the following properties relevant to this paper. Given matrices (or vectors) $\mathbf{U}$, $\mathbf{V}$, $\mathbf{W}$, $\mathbf{X}$, then

1. $(\mathbf{U} \otimes \mathbf{V})^T = \mathbf{U}^T \otimes \mathbf{V}^T$.

2. $(\mathbf{U} \otimes \mathbf{V})(\mathbf{W} \otimes \mathbf{X}) = (\mathbf{U}\mathbf{W}) \otimes (\mathbf{V}\mathbf{X})$

For column vectors $\mathbf{j}$, $\mathbf{k}$, $\mathbf{m}$, $\mathbf{n}$, of dimension d, the previous properties generate the following equalities, which are the bases of the results shown in Sections 4 and 5:

$(\mathbf{j} \otimes \mathbf{k})^T (\mathbf{m} \otimes \mathbf{n}) = (\mathbf{j}^T \otimes \mathbf{k}^T)(\mathbf{m} \otimes \mathbf{n}) =$

$(\mathbf{j}^T\mathbf{m}) \otimes (\mathbf{k}^T\mathbf{n}) = (\mathbf{j}.\mathbf{m}) \otimes (\mathbf{k}.\mathbf{n}) = (\mathbf{j}.\mathbf{m})(\mathbf{k}.\mathbf{n})$.

Scalar products are numbers and the tensor product between numbers is equivalent to a standard arithmetic product.

One way to obtain the results of the action of a logical matrix on uncertain values arises from the general results published by Mizraji (1992). The first step is to define the vector truth values as linear combinations of the vectors $\mathbf{f}$ and $\mathbf{t}$ where the weights are probabilistic scalars:

$\mathbf{u} = \alpha\mathbf{t} + (1-\alpha)\mathbf{f}$ ; $\mathbf{v} = \beta\mathbf{t} + (1-\beta)\mathbf{f}$ ; $\alpha, \beta \in [0,1]$.

We define the projection Pr of an operator on the vector **t**, and since the weights are probabilistic, we have all the information we need. These projections are:

$Pr(\mathbf{N}) = 1 - \alpha$ ,

$Pr(\mathbf{C}) = 1 - \alpha(1-\beta)$ ,

$Pr(\mathbf{K}) = \alpha\beta$ ,

$Pr(\mathbf{A}) = \alpha + \beta - \alpha\beta$ .

From this it follows that for $\alpha, \beta = ½$ the results of Section 5 are obtained.

**Acknowledgment**

I thank Dr. Andrés Pomi our multiple conversations on the topic of this work